\documentclass{amsart}
\usepackage{amsopn}
\usepackage{amssymb, amscd}

\newcommand{\nc}{\newcommand}

\nc{\vg}{\mathfrak{v} } \nc{\wg}{\mathfrak{w} } \nc{\zg}{\mathfrak{z} }
\nc{\ngo}{\mathfrak{n} } \nc{\kg}{\mathfrak{k} } \nc{\mg}{\mathfrak{m} }
\nc{\bg}{\mathfrak{b} } \nc{\ggo}{\mathfrak{g} } \nc{\ggob}{\overline{\mathfrak{g}}
} \nc{\sog}{\mathfrak{so} } \nc{\sug}{\mathfrak{su} } \nc{\spg}{\mathfrak{sp} }
\nc{\slg}{\mathfrak{sl} } \nc{\glg}{\mathfrak{gl} } \nc{\cg}{\mathfrak{c} }
\nc{\rg}{\mathfrak{r} } \nc{\hg}{\mathfrak{h} } \nc{\tg}{\mathfrak{t} }
\nc{\ug}{\mathfrak{u} } \nc{\dg}{\mathfrak{d} } \nc{\ag}{\mathfrak{a} }
\nc{\pg}{\mathfrak{p} } \nc{\sg}{\mathfrak{s} } \nc{\pca}{\mathcal{P}}
\nc{\nca}{\mathcal{N}} \nc{\lca}{\mathcal{L}} \nc{\oca}{\mathcal{O}}
\nc{\mca}{\mathcal{M}} \nc{\tca}{\mathcal{T}} \nc{\aca}{\mathcal{A}}
\nc{\cca}{\mathcal{C}} \nc{\gca}{\mathcal{G}} \nc{\sca}{\mathcal{S}}
\nc{\hca}{\mathcal{H}} \nc{\bca}{\mathcal{B}} \nc{\dca}{\mathcal{D}}
\nc{\val}{\operatorname{val}}

\nc{\vp}{\varphi} \nc{\ddt}{\tfrac{{\rm d}}{{\rm d}t}} \nc{\im}{\mathtt{i}}

\nc{\SO}{\mathrm{SO}} \nc{\Spe}{\mathrm{Sp}} \nc{\Sl}{\mathrm{SL}}
\nc{\SU}{\mathrm{SU}} \nc{\Or}{\mathrm{O}} \nc{\U}{\mathrm{U}} \nc{\Gl}{\mathrm{GL}}
\nc{\Se}{\mathrm{S}} \nc{\Cl}{\mathrm{Cl}} \nc{\Spein}{\mathrm{Spin}}
\nc{\Pin}{\mathrm{Pin}} \nc{\G}{\mathrm{GL}_n(\RR)} \nc{\g}{\mathfrak{gl}_n(\RR)}

\nc{\RR}{{\Bbb R}} \nc{\HH}{{\Bbb H}} \nc{\CC}{{\Bbb C}} \nc{\ZZ}{{\Bbb Z}}
\nc{\FF}{{\Bbb F}} \nc{\NN}{{\Bbb N}} \nc{\QQ}{{\Bbb Q}} \nc{\PP}{{\Bbb P}}

\nc{\vs}{\vspace{.2cm}} \nc{\vsp}{\vspace{1cm}} \nc{\ip}{\langle\cdot,\cdot\rangle}
\nc{\ipp}{(\cdot,\cdot)} \nc{\la}{\langle} \nc{\ra}{\rangle} \nc{\unm}{\tfrac{1}{2}}
\nc{\unc}{\tfrac{1}{4}} \nc{\und}{\tfrac{1}{16}} \nc{\no}{\vs\noindent}
\nc{\lam}{\Lambda^2(\RR^n)^*\otimes\RR^n} \nc{\tangz}{{\rm T}^{\rm Zar}}
\nc{\nor}{{\sf n}}  \nc{\mum}{/\!\!/} \nc{\kir}{/\!\!/\!\!/}
\nc{\Ri}{\tfrac{4\Ric_{\mu}}{||\mu||^2}} \nc{\ds}{\displaystyle}
\nc{\ben}{\begin{enumerate}} \nc{\een}{\end{enumerate}} \nc{\f}{\frac}

\nc{\He}{\operatorname{Hess}} \nc{\ad}{\operatorname{ad}}
\nc{\Ad}{\operatorname{Ad}} \nc{\rank}{\operatorname{rank}}
\nc{\Irr}{\operatorname{Irr}} \nc{\End}{\operatorname{End}}
\nc{\Aut}{\operatorname{Aut}} \nc{\Inn}{\operatorname{Inn}}
\nc{\Der}{\operatorname{Der}} \nc{\Ker}{\operatorname{Ker}}
\nc{\Iso}{\operatorname{I}} \nc{\Diff}{\operatorname{D}} \nc{\Lie}{\operatorname{L}}
\nc{\tr}{\operatorname{tr}} \nc{\dif}{\operatorname{d}}
\nc{\sen}{\operatorname{sen}} \nc{\modu}{\operatorname{mod}}
\nc{\Ric}{\operatorname{R}} \nc{\Ricci}{\operatorname{Ric}}
\nc{\sym}{\operatorname{sym}} \nc{\symac}{\operatorname{sym^{ac}}}
\nc{\symc}{\operatorname{sym^{c}}} \nc{\scalar}{\operatorname{sc}}
\nc{\grad}{\operatorname{grad}} \nc{\ricci}{\operatorname{ric}}
\nc{\ricciac}{\operatorname{ric^{ac}}} \nc{\riccic}{\operatorname{ric^{c}}}
\nc{\riccig}{\operatorname{ric^{\gamma}}} \nc{\Rin}{\operatorname{M}}
\nc{\Le}{\operatorname{L}} \nc{\tang}{\operatorname{T}}
\nc{\level}{\operatorname{level}} \nc{\rad}{\operatorname{r}}
\nc{\abel}{\operatorname{ab}} \nc{\CH}{\operatorname{CH}}
\nc{\mcc}{\operatorname{mcc}} \nc{\Adj}{\operatorname{Adj}}

\theoremstyle{plain}
\newtheorem{theorem}{Theorem}[section]
\newtheorem{proposition}[theorem]{Proposition}

\newtheorem{lemma}[theorem]{Lemma}

\theoremstyle{definition}
\newtheorem{definition}[theorem]{Definition}

\theoremstyle{remark}
\newtheorem{remark}[theorem]{Remark}

\title{Einstein solvmanifolds are standard}

\author{Jorge Lauret}

\address{FaMAF and CIEM, Universidad Nacional de C\'ordoba, C\'ordoba, Argentina}
\email{lauret@famaf.unc.edu.ar}

\thanks{This research was partially supported by grants from CONICET, ANPCyT (Argentina)
and SeCyT (Universidad Nacional de C\'ordoba)}

\begin{document}

\maketitle

\begin{abstract}
We study Einstein manifolds admitting a transitive solvable Lie group of isometries
({\it solvmanifolds}).  It is conjectured that these exhaust the class of noncompact
homogeneous Einstein manifolds. J. Heber \cite{Hbr} has showed that under a simple algebraic condition (he calls such a solvmanifold {\it standard}), Einstein
solvmanifolds have many remarkable structural and uniqueness properties. In this
paper, we prove that any Einstein solvmanifold is standard, by applying a
stratification procedure adapted from one in geometric invariant theory due to F. Kirwan
\cite{Krw1}.
\end{abstract}

\section{Introduction}\label{intro}

The construction of Einstein metrics on manifolds is a classical problem in
differential geometry and general relativity.  A Riemannian manifold is called {\it
Einstein} if its Ricci tensor is a scalar multiple of the metric.  The Einstein
equation $\Ricci(g)=\lambda g$ is a non-linear second order system of partial
differential equations, and a general understanding of the solutions seems far from
being attained (see \cite[11.4]{Brg2}). General existence and non-existence results are hard to obtain, and
it is a natural simplification to impose additional symmetry assumptions, i.e. to
consider metrics admitting a large Lie group of isometries.  In the homogeneous
case, the Einstein equation becomes a subtle system of algebraic equations, and the
following main general question is still open in both, the compact and noncompact cases:

\begin{quote}
Which homogeneous spaces $G/K$ admit a $G$-invariant Einstein Riemannian metric?
\end{quote}

In this paper we shall consider this question in the noncompact case.  We refer to
\cite{BhmWngZll} and the references therein for an account in the compact case.

All the known examples of noncompact homogeneous Einstein manifolds belong to the class of {\it solvmanifolds}, that is, simply connected solvable
Lie groups $S$ endowed with a left invariant metric (see the survey \cite{cruzchica}). According to a long standing
conjecture attributed to D. Alekseevskii (see \cite[7.57]{Bss}), these might exhaust
the class of noncompact homogeneous Einstein manifolds.

On the other hand, all the known examples of Einstein solvmanifolds satisfy the
following additional condition: if $\sg=\ag\oplus\ngo$ is the orthogonal
decomposition of the Lie algebra $\sg$ of $S$ with $\ngo=[\sg,\sg]$, then
$[\ag,\ag]=0$.  A solvmanifold with such a property is called {\it standard}.  For
instance, any solvmanifold of nonpositive sectional curvature is standard (see
\cite{AznWls}).

Standard Einstein solvmanifolds constitute a distinguished class that has been
deeply investigated by J. Heber, who has derived many remarkable structural and
uniqueness results, by assuming only the standard condition (see \cite{Hbr}). We
shall review some of them.  In contrast to the compact case, a standard Einstein
metric is unique up to isometry and scaling among invariant metrics (\cite[Theorem
E]{Hbr}).  Any standard Einstein solvmanifold is isometric to a solvmanifold whose
underlying metric Lie algebra resembles an Iwasawa subalgebra of a semisimple Lie
algebra in the sense that $\ad{A}$ is symmetric and nonzero for any $A\in\ag$,
$A\ne 0$. Moreover, if $H$ denotes the {\it mean curvature vector} of $S$ (i.e.
$\tr{\ad{A}}=\la H,A\ra$ for all $A\in\ag$), then the eigenvalues of
$\ad{H}|_{\ngo}$ form (up to scaling) a set of natural numbers, called
the {\it eigenvalue type} of $S$. There are finitely many such types in each
dimension. Let $\mca$ be the moduli space of all the isometry classes of Einstein
solvmanifolds of a given dimension with scalar curvature equal to $-1$, and let
$\mca_{{\rm st}}$ be the subspace of those which are standard. Then each eigenvalue
type determines a compact pathwise connected component of $\mca_{{\rm st}}$, which
is homeomorphic to a real semialgebraic set. A main result in \cite{Hbr} shows that
$\mca_{st}$ is open in $\mca$ in the $C^{\infty}$-topology (\cite[Theorem G]{Hbr}).

The goal of this paper is to apply an adaptation of a stratification method given in
\cite{Krw1} to prove that actually $\mca_{{\rm st}}=\mca$. In particular, all the
nice structural and uniqueness results in \cite{Hbr} are valid for any
Einstein solvmanifold, and possibly for any noncompact homogeneous Einstein manifold
(if the Alekseevskii's conjecture turns out to be true).

\vspace{.3cm}\noindent {\bf Theorem.} Any Einstein solvmanifold is standard.
\vspace{.3cm}

The proof of the theorem involves a somewhat extensive study of the natural
$\G$-action on the vector space $V_n:=\Lambda^2(\RR^n)^*\otimes\RR^n$, from a geometric
invariant theory point of view (a method already used in \cite[Sections 6.3 and
6.4]{Hbr} in the standard case). We recall that $V_n$ can be viewed as a vector space
containing the space of all $n$-dimensional Lie algebras as an algebraic subset.

We define in Section \ref{st} a $\G$-invariant stratification of $V_n$ satisfying
certain boundary properties (see Theorem \ref{strata}), by adapting a construction
for reductive groups actions on projective algebraic varieties given by F. Kirwan in
\cite[Section 12]{Krw1} in the algebraically closed case (see also \cite{Nss}).  We note that any $\mu\in
V_n$ is {\it unstable} (i.e. $0\in\overline{\G.\mu}$).  The strata are parametrized by
a finite set $\bca$ of diagonal $n\times n$ matrices, and each $\beta\in\bca$ is (up
to conjugation) the `most responsible' direction for the instability of each $\mu$
in the stratum $\sca_{\beta}$, in the sense that $e^{-t\beta}.\mu\to 0$, as
$t\to\infty$ faster than any other one-parameter subgroup having a tangent vector of
the same norm.

We also prove that $\sca_{\beta}$ can be described in terms of semistable (i.e. non
unstable) vectors for a suitable action.  This and the fact that the automorphism
group of any semistable $\mu\in\sca_{\beta}$ must be contained in the parabolic
subgroup of $\G$ defined by $\beta$ are crucial in the proof of several lemmas
needed to prove the theorem.

The first step in the proof of the main theorem, given in Section \ref{eins}, uses
the following fact proved in \cite{minimal}: $S$ is an Einstein solvmanifold with
$\Ricci=cI$ if and only if
\begin{equation}\label{einstein2}
\tr{\left(cI+\unm B+S(\ad{H})\right)E}= \unc\la
\pi(E)[\cdot,\cdot],[\cdot,\cdot]\ra, \qquad\forall E\in\End(\sg),
\end{equation}
where $B$ is the Killing form and $S(\cdot)$ denotes the symmetric part of an operator. Here
$\sg$ is identified with $\RR^m$, the bracket $[\cdot,\cdot]$ of $\sg$ becomes a
vector in $V_m$ and $\pi$ is the representation of $\glg_m(\RR)$ on $V_m$ corresponding to the
$\Gl_m(\RR)$-action.  The next step is to use (\ref{einstein2}) for the right choice of $E$,
namely, $E|_{\ag}=0$, $E|_{\ngo}=\beta+||\beta||^2I$, where $\sca_{\beta}$ is the
stratum the Lie bracket $\mu:=[\cdot,\cdot]|_{\ngo\times\ngo}$ belongs to, as a
vector of $V_n$. Notice that we are identifying $\ngo=[\sg,\sg]$ with $\RR^n$.  We
obtain in this way many expressions which are proved to be nonnegative by using the
lemmas in Section \ref{st}. Finally, by a positivity argument, all these
inequalities turn into equalities, and one of them shows that $[\ag,\ag]=0$.  A key
fact is that $\beta+||\beta||^2I$ is positive definite for any stratum
$\sca_{\beta}$ which meets the closed subset $\nca\subset V_n$ of all nilpotent Lie
brackets.  This is a special feature of this action and is the only point where
we actually use that the vector $\mu=[\cdot,\cdot]|_{\ngo\times\ngo}\in V_n$, to which
we are applying the geometric invariant theory machinery, is a nilpotent Lie algebra.

Partial results on the question if $\mca_{{\rm st}}=\mca$ were obtained by J. Heber
\cite{Hbr}, D. Schueth \cite{Sch}, E. Nikitenko and Y. Nikonorov \cite{NktNkn} and
Y. Nikolayevsky \cite{Nkl}. It is proved for instance in \cite{Nkl} that many
classes of nilpotent Lie algebras can not be the nilradical of a non-standard
Einstein solvmanifold.

More recent results on the structure of standard Einstein solvmanifolds include
interplays with critical points of the square norm of a moment map and Ricci soliton
metrics (see for instance \cite{cruzchica,Nkl3} and the references therein).

We finally mention that the stratification in this paper has also proved to be very
useful in the study of standard Einstein solvmanifolds (see \cite{einsteinsolv}). The
subset $\nca\subset V_n$ parametrizes a set of $(n+1)$-dimensional rank-one (i.e.
$\dim{\ag}=1$) solvmanifolds $\{ S_{\mu}:\mu\in\nca\}$, containing the set of all
those which are Einstein in that dimension. The stratum of $\mu$ determines the
eigenvalue type of a potential Einstein solvmanifold $S_{g.\mu}$, $g\in\G$ (if any),
and so the stratification provides a convenient tool to produce existence results as
well as obstructions for nilpotent Lie algebras to be the nilradical of an Einstein
solvmanifold.

\vs \noindent {\it Acknowledgements.}  The author gratefully acknowledges the many helpful suggestions of Roberto Miatello and Cynthia Will during the preparation of the paper.

\section{A stratification of $V_n$}\label{st}

In this section, we define a $\G$-invariant stratification of a certain real
representation $V_n$ of $\G$ by adapting to our context the construction given by F.
Kirwan in \cite[Section 12]{Krw1} for reductive group representations over an
algebraically closed field.  This construction, in turn, is based on some
instability results due to G. Kempf \cite{Kmp} and W. Hesselink \cite{Hss} (see also \cite{Nss}).  We have
decided to give a self-contained proof of all these results following \cite{Krw1},
bearing in mind that they are crucial in the proof of Theorem \ref{anyst} and that a
direct application of them does not seem feasible.  What is really needed in the proof is the existence of a diagonal matrix $\beta$ satisfying conditions (\ref{adbeta}) and  those stated in Lemmas \ref{betaort}, \ref{delta} and \ref{betapos}.

We consider the vector space
$$
V_n=\lam=\{\mu:\RR^n\times\RR^n\longrightarrow\RR^n : \mu\; \mbox{bilinear and
skew-symmetric}\},
$$
on which there is a natural linear action of $\G$ on the left given by
\begin{equation}\label{action}
g.\mu(X,Y)=g\mu(g^{-1}X,g^{-1}Y), \qquad X,Y\in\RR^n, \quad g\in\G,\quad \mu\in V_n.
\end{equation}

The canonical inner product $\ip$ on $\RR^n$ defines an $\Or(n)$-invariant inner
product on $V_n$ by
\begin{equation}\label{innV}
\la\mu,\lambda\ra=\sum\limits_{ijk}\la\mu(e_i,e_j),e_k\ra\la\lambda(e_i,e_j),e_k\ra,
\end{equation}
where $\{ e_1,...,e_n\}$ is the canonical basis of $\RR^n$.  A Cartan
decomposition for the Lie algebra of $\G$ is given by $\g=\sog(n)\oplus\sym(n)$, that
is, in skew-symmetric and symmetric matrices respectively. We consider the following
$\Ad(\Or(n))$-invariant inner product on $\g$,
\begin{equation}\label{inng}
\la \alpha,\beta\ra=\tr{\alpha \beta^{\mathrm t}}, \qquad \alpha,\beta\in\g.
\end{equation}

\begin{remark} There have been several abuses of notation concerning inner products.  Recall
that $\ip$ has been used to denote an inner product on $\RR^n$, $V_n$ and $\g$.
\end{remark}

The action of $\g$ on $V_n$ obtained by differentiation of (\ref{action}) is given by
\begin{equation}\label{actiong}
\pi(\alpha)\mu=\alpha\mu(\cdot,\cdot)-\mu(\alpha\cdot,\cdot)-\mu(\cdot,\alpha\cdot),
\qquad \alpha\in\g,\quad\mu\in V_n.
\end{equation}

We note that $\pi(\alpha)^t=\pi(\alpha^t)$ for any $\alpha\in\g$.  Let $\tg$ denote
the set of all diagonal $n\times n$ matrices. If $\{ e_1',...,e_n'\}$ is the basis
of $(\RR^n)^*$ dual to the canonical basis then
$$
\{ v_{ijk}=(e_i'\wedge e_j')\otimes e_k : 1\leq i<j\leq n, \; 1\leq k\leq n\}
$$
is a basis of weight vectors of $V_n$ for the action (\ref{action}), where $v_{ijk}$
is actually the bilinear form on $\RR^n$ defined by
$v_{ijk}(e_i,e_j)=-v_{ijk}(e_j,e_i)=e_k$ and zero otherwise.  The corresponding
weights $\alpha_{ij}^k\in\tg$, $i<j$, are given by
$$
\pi(\alpha)v_{ijk}=(a_k-a_i-a_j)v_{ijk}=\la\alpha,\alpha_{ij}^k\ra v_{ijk},
\quad\forall\alpha=\left[\begin{smallmatrix} a_1&&\\ &\ddots&\\ &&a_n
\end{smallmatrix}\right]\in\tg,
$$
where $\alpha_{ij}^k=E_{kk}-E_{ii}-E_{jj}$ and $\ip$ is the inner product defined in
(\ref{inng}).  As usual, $E_{rs}$ denotes the matrix whose only nonzero coefficient
is $1$ at the entry $rs$.  From now on, we will always denote by $\mu_{ij}^k$ the
coefficients of a vector $\mu\in V_n$ with respect to this basis, that is,
$$
\mu=\sum\mu_{ij}^kv_{ijk}, \qquad \mu_{ij}^k\in\RR.
$$
Let $\dca$ denote the set of all $n\times n$ matrices which are diagonalizable, that
is,
$$
\dca=\bigcup_{g\in\G}g\tg g^{-1}.
$$
Consider $m:V_n\times \dca\longrightarrow\RR$ the function defined by letting
$m(\mu,\alpha)$ be the smallest eigenvalue of $\pi(\alpha)$ such that the projection
of $\mu$ onto the corresponding eigenspace is nonzero. Since $\mu_1$ is an
eigenvector of $\pi(\alpha)$ with eigenvalue $a$ if and only if $g.\mu_1$ is an
eigenvector of $\pi(g\alpha g^{-1})$ with eigenvalue $a$, we have that
\begin{equation}\label{strata1}
m(g.\mu,g\alpha g^{-1})=m(\mu,\alpha), \qquad \forall g\in\G.
\end{equation}

It follows from the definition of $m$ that $m(\mu,a\alpha)=am(\mu,\alpha)$ for any
$a>0$ and
$$
m(\mu,\alpha)=\inf\left\{\la\alpha,\alpha_{ij}^k\ra: {\mu_{ij}^k\ne 0}\right\},
\quad\forall\alpha\in\tg.
$$
For each nonzero $\mu\in V_n$ define
$$
Q(\mu)=\inf_{\alpha\in \dca}\{q(\alpha):m(\mu,\alpha)\geq 1\}
$$
and
$$
\Lambda(\mu)=\{\beta\in \dca:q(\beta)=Q(\mu), \quad m(\mu,\beta)\geq 1\},
$$
where $q:\g\longrightarrow\RR$ is defined by $q(\alpha)=\tr{\alpha^2}$.  Note that
$q$ is invariant by $\G$-conjugation, $q(\alpha)>0$ for any nonzero $\alpha\in \dca$
and $q(\alpha)=||\alpha||^2$ for any $\alpha\in\sym(n)$.

\begin{remark}
It would actually be enough to require that $m(\mu,\alpha)=1$ in the definition of $Q(\mu)$, and also it turns out that $m(\mu,\beta)=1$ for any $\beta\in\Lambda(\mu)$, but
the above choice is technically more convenient.
\end{remark}

\begin{remark}
Every $\mu\in V_n$ is unstable for this $\G$-action, i.e. $0\in\overline{\G.\mu}$,
since scalar matrices act as homotheties.  Recall that if $m(\mu,\alpha)>0$ for
$\alpha\in\dca$ then $\lim\limits_{t\to\infty}e^{-t\alpha}.\mu=0$, and thus the
number $Q(\mu)^{-1}$ measures in some sense the degree of instability of $\mu$ as
$$
Q(\mu)^{-\unm}=\sup_{\alpha\in\dca}\{m(\mu,\alpha):q(\alpha)=1\}.
$$
We note that the existence of such a one-parameter subgroup is also necessary for a
vector to be unstable, and is called the numerical criterion of stability. This
notion was established by D. Hilbert to classify homogeneous polynomials, by D.
Mumford \cite{Mmf} in the general algebraically closed case and by R. Richardson and
P. Slodowy \cite{RchSld} in the real case. One therefore may give the following
description of the set $\Lambda(\mu)$: the elements $\beta\in\Lambda(\mu)$ are the
`most responsible' for the instability of $\mu$, in the sense that $e^{-t\beta}.\mu$
converges to zero when $t\to\infty$ more quickly than any other $\alpha$ of the same
norm (recall the definition of $m$).  To show that $\Lambda(\mu)$ lies in a single
conjugacy class (i.e. that such a one-parameter subgroup is essentially unique) will
be actually the main goal of this section. This was proved by G. Kempf \cite{Kmp} in
the complex reductive case.
\end{remark}

\begin{lemma}\label{strata2}
$Q$ is $\G$-invariant, $Q(\mu)>0$ for any $\mu\ne 0$ and
$$
\Lambda(g.\mu)=g\Lambda(\mu)g^{-1}, \qquad \forall\mu\in V_n, \; g\in\G.
$$
\end{lemma}

\begin{proof}
$Q(\mu)$ is always positive for a nonzero $\mu$ since the eigenvalues of
$\pi(\alpha)$ converge to zero when $\alpha$ converges to zero.  It follows from
(\ref{strata1}) that:
$$
\begin{array}{rl}
Q(g.\mu) & =\inf\limits_{\alpha\in \dca}\{ q(\alpha):m(g.\mu,\alpha)\geq 1\} \\ \\
& =\inf\limits_{\alpha\in \dca}\{ q(g^{-1}\alpha g):m(\mu,g^{-1}\alpha
g)\geq 1\}= Q(\mu).
\end{array}
$$
By definition, $\beta\in\Lambda(\mu)$ if and only if $m(\mu,\beta)\geq 1$ and
$q(\beta)\leq q(\alpha)$ for any $\alpha\in \dca$ such that $m(\mu,\alpha)\geq 1$,
which is equivalent to say that $m(g.\mu,g\beta g^{-1})\geq 1$ and $q(g\beta
g^{-1})\leq q(g\alpha g^{-1})$ for any $g\alpha g^{-1}\in \dca$ such that
$m(g.\mu,g\alpha g^{-1})\geq 1$, that is, $g\beta g^{-1}\in\Lambda(g.\mu)$.
\end{proof}

Let $T$ be the subgroup of $\G$ consisting of the diagonal invertible matrices. In
an analogous way, we consider for the action of $T$ on $V_n$ the number
$$
Q_T(\mu)=\inf_{\alpha\in\tg}\{ ||\alpha||^2:m(\mu,\alpha)\geq 1\}
$$
and
$$
\Lambda_T(\mu)=\{\beta\in \tg: ||\beta||^2=Q_T(\mu), \quad m(\mu,\beta)\geq 1\}.
$$

Given a finite subset $X$ of $\tg$, denote by $\CH(X)$ the convex hull of $X$ and by
$\mcc(X)$ the {\it minimal convex combination of} $X$, that is, the (unique) vector
of minimal norm in $\CH(X)$. Each nonzero $\mu\in V_n$ uniquely determines an element
$\beta_{\mu}\in\tg$ given by
$$
\beta_{\mu}=\mcc\left\{\alpha_{ij}^k:\mu_{ij}^k\ne 0\right\}, \qquad
\mu=\sum\mu_{ij}^kv_{ijk}.
$$
We note that $\beta_{\mu}$ is always nonzero since $\tr{\alpha_{ij}^k}=-1$ for all
$i<j$ and consequently $\tr{\beta_{\mu}}=-1$.

\begin{lemma}\label{strata3}
$\Lambda_T(\mu)=\left\{\tfrac{\beta_{\mu}}{||\beta_{\mu}||^2}\right\}$ and
$Q_T(\mu)=\tfrac{1}{||\beta_{\mu}||^2}$.
\end{lemma}

\begin{proof}
Since $\la\beta_{\mu},\alpha_{ij}^k\ra\geq||\beta_{\mu}||^2$ for all $\mu_{ij}^k\ne
0$, we obtain that $m\left(\mu,\tfrac{\beta_{\mu}}{||\beta_{\mu}||^2}\right)\geq 1$.
On the other hand, for $\alpha\in\tg$, if $m(\mu,\alpha)\geq 1$ then
$1\leq\la\alpha,\alpha_{ij}^k\ra$ for all $\mu_{ij}^k\ne 0$ and so
$1\leq\la\alpha,\beta_{\mu}\ra\leq ||\alpha||\; ||\beta_{\mu}||$.  Thus
$||\alpha||\geq\tfrac{1}{||\beta_{\mu}||}=\left|\left|
\tfrac{\beta_{\mu}}{||\beta_{\mu}||^2}\right|\right|$, and the equality holds if and
only if $\alpha=\tfrac{\beta_{\mu}}{||\beta_{\mu}||^2}$.
\end{proof}

\begin{remark}\label{strata4}
Let $T_1$ be any maximal torus of $\G$ and define $\Lambda_{T_1}(\mu)$ as above but
using the Lie algebra of $T_1$ instead of $\tg$.  By considering the weights of the
$T_1$-action on $V_n$ one can prove exactly as above that $\Lambda_{T_1}(\mu)$
consists of a single element: the minimal convex combination of those weights having
nonzero projection of $\mu$ onto their weight spaces.
\end{remark}

If $\mu$ runs through $V_n$, there are only finitely many possible vectors
$\beta_{\mu}$, so we can define for each $\beta\in\tg$ the set
$$
\sca_{\beta}=\Big\{\mu\in V_n\smallsetminus\{ 0\}:\beta\;\mbox{is an element of
maximal norm in}\;\{\beta_{g.\mu}:g\in\G\}\Big\}.
$$
It is clear that $\sca_{\beta}$ is $\G$-invariant for any $\beta\in\tg$,
$$
V_n\smallsetminus\{ 0\}=\bigcup\limits_{\beta\in\tg}\sca_{\beta},
$$
and the set $\{\beta\in\tg:\sca_{\beta}\ne\emptyset\}$ is finite.

\begin{lemma}\label{strata7}
$\tfrac{\beta}{||\beta||^2}\in\Lambda(\mu)$ for all $\mu\in\sca_{\beta}$ such that
$\beta_{\mu}=\beta$. In particular, $Q(\mu)=\tfrac{1}{||\beta||^2}$ for any
$\mu\in\sca_{\beta}$ and
$$
\sca_{\beta}=\G.\left\{\mu\in \sca_{\beta}:\tfrac{\beta}{||\beta||^2}\in\Lambda(\mu)\right\}.
$$
\end{lemma}

\begin{proof}
If $\mu\in V_n$ and $g\in\G$ then from Lemmas \ref{strata2} and \ref{strata3} we
obtain that
$$
\Lambda(\mu)\cap g^{-1}\tg g=g^{-1}\left(\Lambda(g.\mu)\cap\tg\right)g\subset
g^{-1}\Lambda_T(g.\mu)g=\left\{
g^{-1}\tfrac{\beta_{g.\mu}}{||\beta_{g.\mu}||^2}g\right\},
$$
and $m\left(\mu,g^{-1}\tfrac{\beta_{g.\mu}}{||\beta_{g.\mu}||^2}g\right)=m\left(g.\mu,
\tfrac{\beta_{g.\mu}}{||\beta_{g.\mu}||^2}\right)\geq 1$.  Since $\Lambda(\mu)\subset
\dca$ we have that
$$
Q(\mu)=\inf_{g\in\G}\left\{
q\left(g^{-1}\tfrac{\beta_{g.\mu}}{||\beta_{g.\mu}||^2}g\right)\right\}
=\inf_{g\in\G}\left\{\tfrac{1}{||\beta_{g.\mu}||^2}\right\}.
$$
So if $\mu\in\sca_{\beta}$ then $Q(\mu)=\tfrac{1}{||\beta||^2}$ and if in addition
$\beta_{\mu}=\beta$, then $\tfrac{\beta}{||\beta||^2}\in\Lambda(\mu)$ since in this case
$m\left(\mu,\tfrac{\beta}{||\beta||^2}\right)\geq 1$.  The last assertion thus also follows.
\end{proof}

Let us consider the Weyl chamber of $\g$ and its closure, respectively given by
$$
\tg^+=\left\{\left[\begin{smallmatrix} a_1&&\\ &\ddots&\\ &&a_n
\end{smallmatrix}\right]\in\tg:a_1<...<a_n\right\}, \quad \overline{\tg^+}=\left\{\left[\begin{smallmatrix} a_1&&\\ &\ddots&\\ &&a_n
\end{smallmatrix}\right]\in\tg:a_1\leq...\leq a_n\right\}.
$$
For $\alpha\in\overline{\tg^+}$ we define the parabolic subgroup $P_{\alpha}:=B\G_{\alpha}$,
where $B$ is the subgroup of $\G$ of all lower triangular invertible matrices and
$\G_{\alpha}=\{ g\in\G:g\alpha g^{-1}=\alpha\}$.  In general, for $\alpha'\in \dca$,
we let $P_{\alpha'}:=gP_{\alpha}g^{-1}$ if $\alpha'=g\alpha g^{-1}$,
$\alpha\in\overline{\tg^+}$. This is well defined since $h\alpha h^{-1}=g\alpha g^{-1}$ implies
that $h^{-1}g\in\G_{\alpha}\subset P_{\alpha}$ and so
$h^{-1}gP_{\alpha}g^{-1}h=P_{\alpha}$.

It is easy to see that there is an ordered basis of $V_n$ with respect to which the action of $g$ on $V_n$ is
lower triangular for any $g\in B$, and furthermore the eigenvalues of $\pi(\alpha)$
are increasing for any $\alpha\in\overline{\tg^+}$.  It is then easy to check that
$$
m(\mu,g\alpha g^{-1})=m(g.\mu,\alpha)\geq 1 \qquad \forall\alpha\in\overline{\tg^+}, \;
m(\mu,\alpha)\geq 1, \; g\in P_{\alpha},
$$
from which it follows that
\begin{equation}\label{strata5}
    g\alpha g^{-1}\in\Lambda(\mu)\qquad\forall\alpha\in\Lambda(\mu)\cap\overline{\tg^+}, \; g\in P_{\alpha}.
\end{equation}

More generally, if $\alpha'\in\Lambda(\mu)$ then there exists $h\in\G$ such that
$\alpha'=h\alpha h^{-1}$, $\alpha\in\overline{\tg^+}$, and so any $g\in P_{\alpha'}$ is of the
form $g=hg_1h^{-1}$ with $g_1\in P_{\alpha}$.  It follows from Lemma \ref{strata2}
that $\alpha\in\Lambda(h^{-1}.\mu)$ and then again from Lemma \ref{strata2} and
(\ref{strata5}) we get that $g\alpha' g^{-1}=hg_1\alpha
g_1^{-1}h^{-1}\in\Lambda(\mu)$.  Thus
\begin{equation}\label{strata6}
    g\alpha' g^{-1}\in\Lambda(\mu)\qquad\forall\alpha'\in\Lambda(\mu),\; g\in P_{\alpha'}.
\end{equation}

\begin{proposition}\label{pbeta}
For all $\alpha,\beta\in\Lambda(\mu)$ we have that $P_{\alpha}=P_{\beta}$, and any
such $P_{\alpha}$ acts transitively on $\Lambda(\mu)$ by conjugation.
\end{proposition}

\begin{proof} $P_{\alpha}$ and $P_{\beta}$ are parabolic subgroups of $\G$ and so there is a maximal torus $T_1$
of $\G$ contained in $P_{\alpha}\cap P_{\beta}$ (see for instance \cite[Theorem
74.2]{Frd}).  This implies the existence of elements $g\in P_{\alpha}$, $h\in
P_{\beta}$ such that $g\alpha g^{-1}$ and $h\beta h^{-1}$ both lie in Lie($T_1$),
the Lie algebra of $T_1$. It follows from (\ref{strata6}) that $g\alpha g^{-1}$ and
$h\beta h^{-1}$ both belong to $\Lambda(\mu)\cap{\rm
Lie}(T_1)\subset\Lambda_{T_1}(\mu)$, and so $g\alpha g^{-1}=h\beta h^{-1}$ since
$\Lambda_{T_1}(\mu)$ consists of a single element (Remark \ref{strata4}).  Thus
$$
P_{\alpha}=gP_{\alpha} g^{-1}=P_{g\alpha g^{-1}}=P_{h\beta
h^{-1}}=hP_{\beta}h^{-1}=P_{\beta}
$$
and $P_{\alpha}$ acts transitively on $\Lambda(\mu)$.
\end{proof}

\begin{definition}\cite[2.11]{Krw1}
A finite collection $\{ S_i:i\in I\}$ of subsets of a topological space $X$ form a
{\it stratification} of $X$ if $X$ is the disjoint union of the $S_i$, $i\in I$, and
there is a partial order $>$ on the indexing set $I$ such that
$$
\overline{S_i}\subset S_i\cup\bigcup\limits_{j>i}S_j \qquad\forall i\in I.
$$
The $S_i$'s are called the {\it strata} of $X$.
\end{definition}

For each $\beta\in\tg$ we define
$$
W_{\beta}=\{\mu\in
V_n:\la\beta,\alpha_{ij}^k\ra\geq||\beta||^2, \quad\forall\mu_{ij}^k\ne 0\},
$$
that is, the direct sum of all the eigenspaces of $\pi(\beta)$ with eigenvalues
$\geq ||\beta||^2$.  We also consider
$$
\bca=\{\beta\in\overline{\tg^+}:\sca_{\beta}\ne\emptyset\},
$$
and we can now state the main result of this section.

\begin{theorem}\label{strata}
The collection $\{\sca_{\beta}:\beta\in\bca\}$ is a $\G$-invariant stratification of
$V_n\smallsetminus\{ 0\}${\rm :}
\begin{itemize}
\item[(i)] $V_n\smallsetminus\{ 0\}=\bigcup\limits_{\beta\in\bca}\sca_{\beta}$ \quad {\rm (}disjoint union{\rm )}.

\item[(ii)] $\overline{\sca}_{\beta}\smallsetminus\sca_{\beta}\subset
\bigcup\limits_{||\beta'||>||\beta||}S_{\beta'}$, where $\overline{\sca}_{\beta}$ is
the closure of $\sca_{\beta}$ relative to the usual topology of $V_n$.  In particular,
each stratum $\sca_{\beta}$ is a locally closed subset of $V_n\smallsetminus\{ 0\}$.
\end{itemize}
Furthermore, for any $\beta\in\bca$ we have that
\begin{itemize}
\item[(iii)] $W_{\beta}\smallsetminus\{
0\}\subset\sca_{\beta}\cup\bigcup\limits_{||\beta'||>||\beta||}\sca_{\beta'}$.

\item[(iv)] $\sca_{\beta}\cap W_{\beta}=\{\mu\in \sca_{\beta}:\beta_{\mu}=\beta\}$.

\item[(v)] $\sca_{\beta}=\Or(n).\left(\sca_{\beta}\cap W_{\beta}\right)$.
\end{itemize}
\end{theorem}

\begin{proof}
We first prove (i).  Let $\sigma$ be a permutation of $\{ 1,...,n\}$ and define
$g_{\sigma}\in\G$ by $g_{\sigma}e_i=e_{\sigma(i)}$, $i=1,...,n$.  It is easy to see
that if $\mu=\sum\mu_{ij}^kv_{ijk}\in V_n$ then
$$
g_{\sigma}.\mu=\sum\mu_{\sigma^{-1}(i)\sigma^{-1}(j)}^{\sigma^{-1}(k)}v_{ijk},
$$
and since $g_{\sigma}E_{ii}g_{{\sigma}^{-1}}=E_{\sigma(i)\sigma(i)}$ for any $i$ we
have that
$$
\begin{array}{rl}
\{\alpha_{ij}^k:(g_{\sigma}.\mu)_{ij}^k\ne 0\} & =
\{\alpha_{ij}^k:\mu_{\sigma^{-1}(i)\sigma^{-1}(j)}^{\sigma^{-1}(k)}\ne 0\}
=\{\alpha_{\sigma(i)\sigma(j)}^{\sigma(k)}:\mu_{ij}^k\ne 0\} \\ \\
& = \{ g_{\sigma}\alpha_{ij}^kg_{\sigma}^{-1}:\mu_{ij}^k\ne 0\}
=g_{\sigma}\left\{\alpha_{ij}^k:\mu_{ij}^k\ne 0\right\}g_{\sigma}^{-1}.
\end{array}
$$
This implies that
$$
\CH\{\alpha_{ij}^k:(g_{\sigma}.\mu)_{ij}^k\ne 0\}=
g_{\sigma}\CH\left\{\alpha_{ij}^k:\mu_{ij}^k\ne 0\right\}g_{\sigma}^{-1}
$$
and thus $\beta_{g_{\sigma}.\mu}=g_{\sigma}\beta_{\mu}g_{\sigma}^{-1}$ for any
$\mu\in V_n$ and permutation $\sigma$.  We can therefore guarantee the existence of an
element of maximal norm in $\{\beta_{g.\mu}:g\in\G\}$ which lies in $\overline{\tg^+}$ (recall
that $g_{\sigma}\in\Or(n)$ and so $||g_{\sigma}\beta g_{\sigma}^{-1}||=||\beta||$
for any $\beta\in\tg$), and so
$$
V_n\smallsetminus\{ 0\}=\bigcup\limits_{\beta\in\bca}\sca_{\beta}.
$$
Let us now prove the disjointness of this union.  Assume that
$\emptyset\ne\sca_{\beta}\cap\sca_{\beta'}$, $\beta,\beta'\in\bca$. Thus there
exists $\mu\in\sca_{\beta}\cap\sca_{\beta'}$ such that $\beta_{\mu}=\beta$,
$\beta_{g.\mu}=\beta'$ for some $g\in\G$ and $||\beta||=||\beta'||$.  It then
follows from Lemma \ref{strata7} that $\tfrac{\beta}{||\beta||^2}\in\Lambda(\mu)$
and $\tfrac{\beta'}{||\beta'||^2}\in\Lambda(g.\mu)$, or equivalently,
$g^{-1}\tfrac{\beta'}{||\beta'||^2}g\in\Lambda(\mu)$ (Lemma \ref{strata2}).  By
Proposition \ref{pbeta} we get that $\beta$ and $\beta'$ are conjugate and therefore
$\beta=\beta'$ since both are in $\overline{\tg^+}$.

We now prove part (iii).  For any $\mu\in W_{\beta}$ we have that
$\la\beta,\beta_{\mu}\ra\geq||\beta||^2$, and hence either $\beta=\beta_{\mu}$ or
$||\beta||<||\beta_{\mu}||$.  Thus if $\mu\in\sca_{\beta_{\mu}}$ we are done.
Otherwise, $\mu\in\sca_{\beta'}$ for some $\beta'\in\overline{\tg^+}$ such that $||\beta'||\geq
||\beta_{\mu}||$, and so $||\beta'||>||\beta||$ and (iii) follows.

To prove (ii) recall first that for each $\mu\in\sca_{\beta}$ there exists $g\in G$
such that $\beta_{g.\mu}=\beta$ and thus $g.\mu\in W_{\beta}$.  This implies that
$\sca_{\beta}\subset \G W_{\beta}$, a closed subset since $W_{\beta}$ is a subspace
and $\G/P_{\beta}$ is compact (recall that $W_{\beta}$ is $P_{\beta}$-invariant),
and hence $\overline{\sca}_{\beta}\subset\G W_{\beta}$.  It now follows from (iii)
 that
$$
\overline{\sca}_{\beta}\smallsetminus\sca_{\beta}\subset\G
W_{\beta}\smallsetminus\sca_{\beta}\subset\bigcup_{||\beta'||>||\beta||}\sca_{\beta'},
$$
as was to be shown.  The last assertion in (ii) follows from:
$$
\sca_{\beta}=\overline{\sca}_{\beta}\smallsetminus\bigcup\limits_{||\beta'||>||\beta||}S_{\beta'}
$$
and $\bigcup\limits_{||\beta'||>||\beta||}S_{\beta'}$ is also a closed subset of
$V_n\smallsetminus\{ 0\}$.  Thus $\sca_{\beta}$ is the intersection of an open subset
and a closed subset, that is, $\sca_{\beta}$ is locally closed.

We now prove (iv). If $\mu\in W_{\beta}$ then $(\beta,\beta_{\mu})\geq ||\beta||^2$,
and if in addition $\mu\in\sca_{\beta}$, then $||\beta_{\mu}||\leq ||\beta||$.  Thus
$\beta=\beta_{\mu}$ for any $\mu\in W_{\beta}\cap \sca_{\beta}$. Conversely, if
$\mu\in\sca_{\beta}$ and $\beta_{\mu}=\beta$, then $\mu\in W_{\beta}$ since
$\beta_{\mu}=\mcc\{\alpha_{ij}^k:\mu_{ij}^k\ne 0\}$.

Finally, (v) follows from the fact that $W_{\beta}$ is $B$-invariant and hence
$$
\sca_{\beta}=\G.(\sca_{\beta}\cap W_{\beta})=(\Or(n)B).(\sca_{\beta}\cap
W_{\beta})=\Or(n)(\sca_{\beta}\cap W_{\beta}).
$$
This concludes the proof of the theorem.
\end{proof}

It follows from Proposition \ref{pbeta} that for each nonzero $\mu\in V_n$ there
exists a parabolic subgroup $P_\mu\subset \G$ acting transitively on $\Lambda(\mu)$,
which satisfies $P_{\mu}=P_{\alpha}$ for any $\alpha\in\Lambda(\mu)$.  If
$\alpha\in\Lambda(\mu)$ and $g\in\G$ satisfies $g\alpha g^{-1}\in\Lambda(\mu)$, then
there exists $h\in P_{\mu}$ such that $hg\alpha g^{-1}h^{-1}=\alpha$.  Thus
$hg\in\G_{\alpha}\subset P_{\mu}$ and so $g\in P_{\mu}$.  This implies that
$$
P_\mu=\{ g\in \G:\Ad(g)\alpha\in\Lambda(\mu)\}, \qquad\forall\alpha\in\Lambda(\mu),
$$
which in turn gives $\Aut(\mu)\subset P_{\mu}$, where $\Aut(\mu)$ is the
automorphism group of the algebra $\mu$.  Indeed, $\Lambda(\mu)$ is
$\Aut(\mu)$-invariant since $m(\mu,\alpha)=m(g.\mu,g\alpha g^{-1})=m(\mu,g\alpha
g^{-1})$ for all $g\in\Aut(\mu)$.  We therefore obtain that
\begin{equation}\label{isocont}
\Der(\mu)\subset \pg_{\mu},
\end{equation}
where $\Der(\mu)=\{\alpha\in\g:\pi(\alpha)\mu=0\}$ is the Lie algebra of
derivations of $\mu$ and $\pg_{\mu}$ is the Lie algebra of $P_{\mu}$.  We note that
if $\mu\in\sca_{\beta}\cap W_\beta$, $\beta\in\bca$, then
$\tfrac{\beta}{||\beta||^2}\in\Lambda(\mu)$ (see Lemma \ref{strata7}) and so
$P_{\mu}=P_{\beta}=B\G_{\beta}$. It is then easy to check by using (\ref{isocont})
that
\begin{equation}\label{adbeta}
\la[\beta,D],D\ra\geq 0, \qquad \forall D\in\Der(\mu),\quad\mu\in\sca_{\beta}\cap W_\beta.
\end{equation}

We will now give a description of the strata in terms of semistable vectors.  For
each $\beta\in\tg$ consider the sets
$$
\begin{array}{l}
Z_{\beta}=\{ \mu\in V_n:\la\beta,\alpha_{ij}^k\ra=||\beta||^2,\quad\forall \mu_{ij}^k\ne 0\}, \\ \\
Y_{\beta}=\{ \mu\in W_{\beta}:\la\beta,\alpha_{ij}^k\ra=||\beta||^2\;\mbox{for at
least one}\; \mu_{ij}^k\ne 0\}.
\end{array}
$$
Thus $Z_{\beta}\subset Y_{\beta}\subset W_{\beta}$, and $Z_{\beta}$ is actually the
eigenspace of $\pi(\beta)$ with eigenvalue $||\beta||^2$.  $Z_{\beta}$ is therefore
$\G_{\beta}$-invariant and since $W_{\beta}$ is so, $Y_{\beta}$ turns to be
$\G_{\beta}$-invariant as well.  Let $\ggo_{\beta}$ denote the Lie algebra of
$\G_{\beta}$, that is, $\ggo_{\beta}=\{\alpha'\in\g:[\alpha',\alpha]=0\}$.

\begin{lemma}\label{zbeta}
For any $\mu\in Z_{\beta}$, $\Lambda(\mu)\cap\ggo_{\beta}\ne\emptyset$.  In
particular, there exists $g\in \G_{\beta}$ such that
$Q(\mu)=\tfrac{1}{||\beta_{g.\mu}||^2}$.
\end{lemma}

\begin{proof}
If $\mu\in Z_{\beta}$ then $\beta+||\beta||^2I\in\Der(\mu)$ and so
$\beta+||\beta||^2I\in\pg_{\mu}$ (see (\ref{isocont})), where $\Lambda(\mu)$ is
contained. Thus for any $\alpha\in\Lambda(\mu)$ there exists $h\in P_\mu$ such that
$[h\alpha h^{-1},\beta]=0$, and so $h\alpha h^{-1}\in\Lambda(\mu)\cap\ggo_{\beta}$
and the first assertion follows. For the second one, we first note that if
$\gamma\in\Lambda(\mu)\cap\ggo_{\beta}$ then there exists $g\in \G_{\beta}$ such
that $g\gamma g^{-1}\in\Lambda(\mu)\cap\tg\subset\Lambda_T(\mu)$, which by Lemma
\ref{strata3} gives $g\gamma g^{-1}=\tfrac{\beta_{g.w}}{||\beta_{g.w}||^2}$.  Thus
$$
Q(\mu)=q(\gamma)=q(g\gamma g^{-1})=\tfrac{1}{||\beta_{g.\mu}||^2},
$$
as asserted.
\end{proof}

\begin{proposition}\label{zybeta}
For any $\beta\in\bca$ the $\G_{\beta}$-invariant subsets
$Z_{\beta}^{ss}:=Z_{\beta}\cap\sca_{\beta}$ and
$Y_{\beta}^{ss}:=Y_{\beta}\cap\sca_{\beta}$ satisfy:
\begin{itemize}
\item[(i)] $Y_{\beta}^{ss}=\sca_\beta\cap W_\beta$; in particular $\sca_{\beta}=\Or(n).Y_{\beta}^{ss}$.

\item[(ii)] $Y_{\beta}^{ss}=\{ \mu\in Y_{\beta}:p_{\beta}(\mu)\in Z_{\beta}^{ss}\}$, where
$p_{\beta}:W_{\beta}\longrightarrow Z_{\beta}$ is the orthogonal projection on
$Z_{\beta}$.
\end{itemize}
\end{proposition}

\begin{proof}
It is clear that $Y_{\beta}^{ss}\subset\sca_{\beta}\cap W_\beta$ as $Y_{\beta}\subset W_\beta$. Conversely, if $\mu\in\sca_{\beta}\cap W_\beta$, then $\beta_\mu=\beta$ and it is not possible to have
$\la\beta,\alpha_{ij}^k\ra>||\beta||^2$ for any $\mu_{ij}^k\ne 0$ since this would
contradict the fact that $\beta=\mcc\{\alpha_{ij}^k:\mu_{ij}^k\ne 0\}$.  This
implies that $\mu\in Y_{\beta}$ and so in $Y_{\beta}^{ss}$, and hence the first assertion in (i) follows.  The second one follows from Theorem \ref{strata}, (v).

To prove (ii), we first note that if $\mu\in Y_{\beta}$ and $\mu\notin\sca_{\beta}$
then $\mu$ must lie in one of the strata $\sca_{\beta'}$ with $||\beta'||>||\beta||$
as $\mu\in W_{\beta}$ (see Theorem \ref{strata}, (iii)). We have that
\begin{equation}\label{limpbeta}
\lambda:=p_{\beta}(\mu)=\lim_{t\to\infty}e^{-t(\beta+||\beta||^2I)}.\mu\in\overline{\G.\mu},
\end{equation}
but then $\lambda\in\overline{\sca}_{\beta'}$ and so
$\lambda\notin\sca_{\beta}$ by Theorem \ref{strata}, (ii). The set on the right hand
side is then contained in $Y_{\beta}^{ss}$. Conversely, if $\mu\in Y_{\beta}^{ss}$
and $\lambda=p_{\beta}(\mu)$ then it follows from Theorem \ref{strata}, (iv) that
$\beta=\beta_\mu=\beta_\lambda$. Let us assume that $\lambda\in\sca_{\beta'}$ with
$||\beta'||>||\beta||$ and $\beta'=\beta_{h.\lambda}$ for some $h\in \G$.  From
Lemma \ref{zbeta} we obtain that $Q(\lambda)^{-1}=||\beta_{g.\lambda}||^2$ for some
$g\in \G_{\beta}$, and so $g.\lambda\in Z_{\beta}$ and
$\beta_{g.\lambda}=\beta_{g.\mu}$, which gives
$$
||\beta||^2\geq
||\beta_{g.\mu}||^2=||\beta_{g.\lambda}||^2=\tfrac{1}{Q(\lambda)}\geq
||\beta_{h.\lambda}||^2=||\beta'||^2,
$$
a contradiction.  Thus, since $\lambda\in\overline{\sca_{\beta}}$, it follows from
Theorem \ref{strata}, (ii) that $\lambda\in\sca_{\beta}$, as claimed.
\end{proof}

Let $H_{\beta}$ be the connected Lie subgroup of $\G$ with Lie algebra
$\hg_{\beta}$, the orthogonal complement of $\beta$ in $\ggo_{\beta}$.  If $\beta=\beta_\mu$ for some nonzero $\mu\in V_n$, then all its entries are in $\QQ$, as $\beta$ is the minimal convex combination of a subset of $\{\alpha_{ij}^k\}$ and all these matrices have rational entries. In particular, $\beta$ is rational and $H_\beta$ is therefore a real reductive algebraic group for any $\beta\in\tg$ such that $\sca_\beta\ne\emptyset$.  In that case,
$\hg_{\beta}=(\sog(n)\cap\ggo_{\beta})\oplus(\sym(n)\cap\hg_{\beta})$ is a Cartan
decomposition and $\tg\cap\hg_{\beta}=\{\alpha\in\tg:\la\alpha,\beta\ra=0\}$ is a
maximal abelian subalgebra of $\sym(n)\cap\hg_{\beta}$.

\begin{definition}
A vector $\mu\in V_n$ is called $H_{\beta}$-{\it semistable} if
$0\notin\overline{H_{\beta}.\mu}$.
\end{definition}

\begin{proposition}\label{zybetass}
For any $\beta\in\tg$ such that $\sca_\beta\ne\emptyset$, the following conditions hold:
\begin{itemize}
\item[(i)] $Z_{\beta}^{ss}$ is the set of $H_{\beta}$-semistable vectors in $Z_{\beta}$.

\item[(ii)] $Y_{\beta}^{ss}$ is the set of $H_{\beta}$-semistable vectors in $W_{\beta}$.
\end{itemize}
\end{proposition}

\begin{proof}
We first prove (i).  Assume that a $\mu\in Z_{\beta}^{ss}$ is not
$H_{\beta}$-semistable.  Thus $\{ 0\}$ is the only closed orbit in
$\overline{H_{\beta}.\mu}$ (see \cite[9.3]{RchSld}) and so there exists
$\alpha\in\sym(n)\cap\hg_{\beta}$ such that
$\lim\limits_{t\to\infty}e^{-t\alpha}.\mu=0$ (see \cite[Lemma 3.3]{RchSld}).  We
take $g\in\Or(n)_{\beta}$ such that $g\alpha g^{-1}\in\tg$, but then $g.\mu\in
Z_{\beta}^{ss}$ and $\lim\limits_{t\to\infty}e^{-tg\alpha g^{-1}}.(g.\mu)=0$ as
well, from which it follows that we can assume $\alpha\in\tg$.  This implies that
$\lim\limits_{t\to\infty}\sum\mu_{ij}^ke^{-t\la\alpha,\alpha_{ij}^k\ra}v_{ijk}=0$
and consequently $\la\alpha,\alpha_{ij}^k\ra>0$ for all $\mu_{ij}^k\ne 0$.  Thus
$\la\alpha,\beta_{\mu}\ra>0$, a contradiction, since $\beta=\beta_{\mu}$ by Theorem
\ref{strata}, (iv) and hence $\alpha\perp\beta$.

Conversely, if $\mu\in Z_{\beta}$ is $H_{\beta}$-semistable and
$\mu\notin\sca_{\beta}$ then $\mu\in\sca_{\beta'}$ with $||\beta'||>||\beta||$
(recall that $\mu\in W_{\beta}$ and see Theorem \ref{strata}, (iii)).  If $g\in
\G_{\beta}$ and $\beta_{g.\mu}=\alpha+a\beta$ for a nonzero $\alpha\in\tg$,
$\alpha\perp\beta$, then
$$
0=\lim_{t\to\infty}\sum
(g.\mu)_{ij}^ke^{-t\la\beta_{g.\mu}-a\beta,\alpha_{ij}^k\ra}v_{ijk}=
\lim_{t\to\infty}e^{-t\alpha}.(g.\mu)\in\overline{H_{\beta}.\mu},
$$
since $g.\mu\in Z_{\beta}\cap\sca_{\beta'}$ and hence
$$
\la\beta_{g.\mu},\alpha_{ij}^k\ra\geq ||\beta_{g.\mu}||^2> ||a\beta||^2=\la
a\beta,\alpha_{ij}^k\ra, \qquad\forall (g.\mu)_{ij}^k\ne 0.
$$
This contadicts the fact that $\mu$ is $H_{\beta}$-semistable and so $\beta_{g.\mu}$
is a scalar multiple of $\beta$.  This gives $\beta_{g.\mu}=\beta$ for any $g\in
\G_{\beta}$ as $\tr{\beta_{g.\mu}}=\tr{\beta}=-1$.  It then follows from Lemma
\ref{zbeta} that $Q(\mu)^{-1}=||\beta||^2$, a contradiction, since
$Q(\mu)^{-1}=||\beta'||^2$ by Lemma \ref{strata7}.  This concludes the proof of (i).

Let $\mu\in W_{\beta}$.  In order to prove (ii), we must show that
$\mu\in\sca_{\beta}$ if and only if $\mu$ is $H_{\beta}$-semistable, and since both
conditions imply as in (i) that $\beta_\mu=\beta$, we may assume that $\mu\in
Y_{\beta}$.  If $\lambda=p_{\beta}(\mu)\in Z_{\beta}$ then
$\beta_\lambda=\beta_\mu$, and since $p_{\beta}(g.\mu)=g.p_{\beta}(\mu)$ for any
$g\in \G_{\beta}$ we obtain that $\mu$ is $H_{\beta}$-semistable if and only if
$\lambda$ is so.  Thus (ii) follows from (i) and the fact that $\mu\in\sca_{\beta}$
if and only if $\lambda\in\sca_{\beta}$ (see Proposition \ref{zybeta}, (ii)).
\end{proof}

In what follows, we prove a series of lemmas which will be needed in the proof of
Theorem \ref{anyst}.

\begin{lemma}\label{betaort}
If $\mu\in Y_{\beta}^{ss}$ then $\tr{\beta D}=0$ for any $D\in\Der(\mu)$.
\end{lemma}

\begin{proof}
It follows from (\ref{isocont}) that $D\in\pg_{\mu}=\bg+\g_{\beta}$, where $\bg$ is
the Lie subalgebra of $\g$ of all lower triangular matrices.  Thus
$\lim\limits_{t\to\infty}e^{-t\beta}De^{t\beta}$ exists and it is an element
$A\in\g_{\beta}$.  If $\lambda=p_{\beta}(\mu)$ then we can use (\ref{limpbeta}) to
show that
$$
\pi(A)\lambda=\lim_{t\to\infty}\pi(e^{-t\beta}De^{t\beta})e^{-t(\beta+||\beta||^2I)}.\mu=
\lim_{t\to\infty}e^{t||\beta||^2}e^{-t\beta}.\pi(D)\mu=0,
$$
that is, $A\in\Der(\lambda)$.  We decompose $A$ as $A=a\beta+A'$ with
$A'\perp\beta$. By using that $e^{tD}.\mu=\mu$ for all $t$ we obtain that
$e^{tA'}.\lambda=e^{-ta\beta}.\lambda=e^{-ta||\beta||^2}\lambda$ (recall that
$\lambda\in Z_{\beta}^{ss}$ by Proposition \ref{zybeta}, (ii)) and hence $a=0$ since
otherwise $0\in\overline{H_{\beta}.\lambda}$, contradicting the fact that $\lambda$
is $H_{\beta}$-semistable (see Proposition \ref{zybetass}, (ii)). This implies that
$\tr{\beta D}=\tr{\beta A}=\tr{\beta A'}=0$, as was to be shown.
\end{proof}

\begin{lemma}\label{delta}
$\la\pi(\beta+||\beta||^2I)\mu,\mu\ra\geq 0$ for any $\mu\in W_{\beta}$.
\end{lemma}

\begin{proof}
If $\mu=\sum\mu_{ij}^kv_{ij}^k\in W_{\beta}$ then $\la\beta,\alpha_{ij}^k\ra\geq
||\beta||^2$ for all $\mu_{ij}^k\ne 0$ and henceforth
$$
\la\pi(\beta+||\beta||^2I)\mu,\mu\ra=\la\pi(\beta)\mu,\mu\ra-||\beta||^2||\mu||^2
=\sum(\mu_{ij}^k)^2\la\beta,\alpha_{ij}^k\ra-||\beta||^2||\mu||^2\geq 0,
$$
as claimed.
\end{proof}

The space of all $n$-dimensional nilpotent Lie algebras can be parametrized by the set
$$
\nca=\{\mu\in V_n:\mu\;\mbox{satisfies the Jacobi identity and is nilpotent}\},
$$
which is an algebraic subset of $V_n$ as the Jacobi identity and the nilpotency
condition can both be expressed as zeroes of polynomial functions.  Note that $\nca$
is $\G$-invariant and Lie algebra isomorphism classes are precisely $\G$-orbits.

\begin{lemma}\label{betapos}
$\beta+||\beta||^2$ is positive definite for every $\beta\in\bca$ such that
$\sca_{\beta}\cap\nca\ne\emptyset$.
\end{lemma}

\begin{proof}
By Proposition \ref{zybeta}, (i) there exists $\mu\in Y_{\beta}^{ss}\cap\nca$, which
is therefore $H_{\beta}$-semistable by Proposition \ref{zybetass}, (ii).  Thus there
exists a nonzero $\lambda\in\overline{H_{\beta}.\mu}\subset Y_{\beta}^{ss}\cap\nca$
such that $||\lambda||\leq ||\lambda'||$ for any
$\lambda'\in\overline{H_{\beta}.\mu}$.  Let us consider $\Ric_{\lambda}\in\g$
defined implicitly by
$$
\la\Ric_{\lambda},\alpha\ra=\tfrac{1}{4}\la\pi(\alpha)\lambda,\lambda\ra, \qquad
\forall\alpha\in\g.
$$
We note that $\Ric_{\lambda}\in\sym(n)$ as for any $\alpha\in\sog(n)$,
$\la\pi(\alpha)\lambda,\lambda\ra=\unm\ddt|_{t=0}||e^{t\alpha}.\lambda||^2=0$.  Since
$\lambda$ is a vector of minimal norm in $H_{\beta}.\lambda$ we obtain
$$
\la\Ric_{\lambda},\alpha\ra=\unc\ddt|_{t=0}\la
e^{t\alpha}.\lambda,\lambda\ra=\tfrac{1}{4}\ddt|_{t=0}||e^{\tfrac{t}{2}\alpha}.\lambda||^2=0,
\qquad\forall\alpha\in\hg_{\beta}\cap\sym(n),
$$
and hence the orthogonal projection of $\Ric_{\lambda}$ on $\ggo_{\beta}$ is a
scalar multiple of $\beta$. Thus, such a projection equals
$\tfrac{||\lambda||^2}{4}\beta$ as $\tr{\Ric_{\lambda}}=\la\Ric_{\lambda},I\ra
=\unc\la\pi(I)\lambda,\lambda\ra=-\tfrac{||\lambda||^2}{4}$ (recall that
$\tr{\beta}=-1$).  This implies that if $\lambda=\sum\lambda_{ij}^kv_{ijk}$ then
\begin{equation}\label{betapos1}
||\beta||^2=\tfrac{4}{||\lambda||^2}\la\Ric_{\lambda},\beta\ra
=\tfrac{1}{||\lambda||^2}\la\pi(\beta)\lambda,\lambda\ra=
\sum\tfrac{(\lambda_{ij}^k)^2}{||\lambda||^2}\la\beta,\alpha_{ij}^k\ra.
\end{equation}

But $\lambda\in Y_{\beta}^{ss}\subset W_{\beta}$ and so
$$
\la\beta,\alpha_{ij}^k\ra\geq ||\beta||^2, \qquad\forall\lambda_{ij}^k\ne 0,
$$
which implies by (\ref{betapos1}) that these must all be equalities and hence
$\lambda\in Z_{\beta}$.

On the other hand, for any $\alpha\in\g$ we have that
$$
\la[\Ric_{\lambda},\beta],\alpha\ra =-\la\Ric_{\lambda},[\beta,\alpha]\ra
=-\unc\la\pi([\beta,\alpha])\lambda,\lambda\ra=0
$$
since $\pi(\beta)$ is symmetric and $\pi(\beta)\lambda=||\beta||^2\lambda$.  We
therefore obtain that $\Ric_{\lambda}\in\ggo_{\beta}$ and so
$\Ric_{\lambda}=\tfrac{||\lambda||^2}{4}\beta$.

If $D:=\tfrac{||\lambda||^2}{4}(\beta+||\beta||^2I)$ and $DX=dX$,
$X\in\RR^n\smallsetminus\{ 0\}$, then since $D\in\Der(\lambda)$,
$$
d\ad_{\lambda}{X}=\ad_{\lambda}(DX)=[D,\ad_{\lambda}{X}]=[\Ric_{\lambda},\ad_{\lambda}{X}],
$$
where $\ad_{\lambda}$ denotes the adjoint representation of the Lie algebra
$\lambda$.  Thus
$$
\begin{array}{rl}
d\tr{\ad_{\lambda}{X}(\ad_{\lambda}{X})^t}
&=\tr{[\Ric_{\lambda},\ad_{\lambda}{X}](\ad_{\lambda}{X})^t}
=\la\Ric_{\lambda},[\ad_{\lambda}{X},(\ad_{\lambda}{X})^t]\ra \\ \\
&=\unc\la\pi([\ad_{\lambda}{X},\ad_{\lambda}{X})^t])\lambda,\lambda\ra
=\unc||\pi((\ad_{\lambda}{X})^t)\lambda||^2,
\end{array}
$$
which implies that $d\geq 0$ as long as $\ad_{\lambda}{X}\ne 0$.  If
$\ad_{\lambda}{X}=0$, since
$\la\Ric_{\lambda}X,X\ra=(d-\tfrac{||\beta||^2||\lambda||^2}{4})||X||^2$ and
\begin{equation}\label{R2}
\la\Ric_{\lambda}X,X\ra=
-\unm\displaystyle{\sum\limits_{ij}}\la\lambda(X,e_i),e_j\ra^2
+\unc\displaystyle{\sum\limits_{ij}}\la\lambda(e_i,e_j),X\ra^2,
\end{equation}
(see \cite[Propositions 3.5, 4.2]{minimal}) we obtain that $d>0$.  If
$\ad_{\lambda}{X}\ne 0$ and $d=0$ then $(\ad_{\lambda}{X})^t\in\Der(\lambda)$, which
is a contradiction since $\lambda$ is a nilpotent Lie algebra (consider the
orthogonal decomposition $\RR^n=\ngo_1\oplus...\oplus\ngo_r$ such that
$\ngo_s\oplus...\oplus\ngo_r$, $s=1,...,r$ is the central descendent series).  We have therefore obtained that in any case $d>0$, and hence $\beta+||\beta||^2I$ is positive definite.
\end{proof}

\begin{remark} Notice that Lemma \ref{betapos} is the only result in this section where we need $\mu$ to be
a nilpotent Lie algebra, and not just any vector in $V_n$.  It is known for instance that semisimple Lie algebras lie
in the stratum $\sca_{\beta}$ for $\beta=-\tfrac{1}{n}I$, and consequently
$\beta+||\beta||^2I=0$ (see \cite{strata}).
\end{remark}

\section{Einstein solvmanifolds}\label{eins}

We now apply the results obtained in Section \ref{st} to prove our main theorem,
namely that Einstein solvmanifolds are all standard.

Let $S$ be a solvmanifold, that is, a simply connected solvable Lie group endowed
with a left invariant Riemannian metric.  Let $\sg$ be the Lie algebra of $S$ and
let $\ip$ denote the inner product on $\sg$ determined by the metric.  We consider
the orthogonal decomposition $\sg=\ag\oplus\ngo$, where $\ngo=[\sg,\sg]$. A
solvmanifold $S$ is called {\it standard} if $[\ag,\ag]=0$.

The mean curvature vector of $S$ is the only element $H\in\ag$ which satisfies
$\la H,A\ra=\tr{\ad{A}}$ for any $A\in\ag$.  If $B$ denotes the symmetric map
defined by the Killing form of $\sg$ relative to $\ip$ then $B(\ag)\subset\ag$ and
$B|_{\ngo}=0$ as $\ngo$ is contained in the nilradical of $\sg$.  The Ricci operator
$\Ricci$ of $S$ is given by (see for instance \cite[7.38]{Bss}):
\begin{equation}\label{ricci}
\Ricci=R-\unm B-S(\ad{H}),
\end{equation}
where $S(\ad{H})=\unm(\ad{H}+(\ad{H})^t)$ is the symmetric part of $\ad{H}$
and $R$ is the symmetric operator defined by
\begin{equation}\label{R}
\la Rx,y\ra=-\unm\displaystyle{\sum\limits_{ij}}\la [x,x_i],x_j\ra\la [y,x_i],x_j\ra
+\unc\displaystyle{\sum\limits_{ij}}\la [x_i,x_j],x\ra\la [x_i,x_j],y\ra,
\end{equation}
for all $x,y\in\sg$, where $\{ x_i\}$ is any orthonormal basis of
$(\sg,\ip)$.  A solvmanifold $S$ is called {\it Einstein} if $\Ricci=cI$ for some
$c\in\RR$.  We refer to \cite{Bss} for a detailed exposition on Einstein manifolds
(see also the surveys in \cite{LbrWng} and \cite[11.4]{Brg2}).

It is proved in \cite[Propositions 3.5, 4.2]{minimal} that $R$ is the only symmetric
operator on $\sg$ such that
\begin{equation}\label{Rmm}
\tr{RE}=\unc\la \pi(E)[\cdot,\cdot],[\cdot,\cdot]\ra, \qquad\forall E\in\End(\sg),
\end{equation}
where we are identifying $\sg$ with $\RR^m$, $[\cdot,\cdot]$ with a vector
in $V_m$, and so $\ip$ is the inner product defined in (\ref{innV}) and $\pi$ is the
representation given in (\ref{actiong}) (see the notation in Section \ref{st}).

We therefore obtain from (\ref{ricci}) and (\ref{Rmm}) that $S$ is an Einstein
solvmanifold with $\Ricci=cI$, if and only if, for any $E\in\End(\sg)$,
\begin{equation}\label{einstein}
\tr{\left(cI+\unm B+S(\ad{H})\right)E}= \unc\sum_{ij}\la
E[x_i,x_j]-[Ex_i,x_j]-[x_i,Ex_j],[x_i,x_j]\ra.
\end{equation}

We are now in a position to prove the main result of this paper.

\begin{theorem}\label{anyst}
Any Einstein solvmanifold is standard.
\end{theorem}

\begin{proof}
Let $S$ be an Einstein solvmanifold with $\Ricci=cI$.  We can assume that $S$ is not
unimodular by using the following fact proved by I. Dotti \cite{Dtt}: a unimodular
Einstein solvmanifold must be flat and consequently standard (see \cite[Proposition
4.9]{Hbr}).  Thus $H\ne 0$ and $\tr{\ad{H}}=||H||^2>0$.  By letting $E=\ad{H}$ in
(\ref{einstein}) and using that $\ad{H}\in\Der(\sg)$ we get
\begin{equation}\label{c}
c=-\tfrac{\tr{S(\ad{H})^2}}{\tr{S(\ad{H})}}.
\end{equation}

In order to apply the results in Section \ref{st}, we identify $\ngo$ with $\RR^n$
via an orthonormal basis $\{ e_1,...,e_n\}$ of $\ngo$ and we set
$\mu:=[\cdot,\cdot]|_{\ngo\times\ngo}$.  In this way, $\mu$ can be viewed as an
element of $\nca\subset V_n$. If $\mu\ne 0$ then $\mu$ lies in a unique stratum
$\sca_{\beta}$, $\beta\in\bca$, by Theorem \ref{strata}, (i).  It then follows from
Proposition \ref{zybeta}, (i) that there exists $g\in\Or(n)$ such that $g.\mu\in
Y_{\beta}^{ss}$. Let $\tilde{g}$ denote the orthogonal map of $(\sg,\ip)$ defined by
$\tilde{g}|_{\ag}=I$, $\tilde{g}|_{\ngo}=g$. We let $\tilde{S}$ to be the
solvmanifold whose Lie algebra $\tilde{\sg}$ is $\sg$ as a vector space and has Lie
bracket
$$
\tilde{g}.[\cdot,\cdot]=\tilde{g}[\tilde{g}^{-1}\cdot,\tilde{g}^{-1}\cdot].
$$
The left invariant metric on $\tilde{S}$ is determined by the same inner product
$\ip$ on $\sg$. We therefore have that $\tilde{S}$ is isometric to $S$, as
$\tilde{g}$ is an isometric isomorphism between the two metric Lie algebras. Thus
$\tilde{S}$ is also Einstein, and since $S$ is standard if and only if $\tilde{S}$
is so, we can assume that $\mu\in Y_{\beta}^{ss}=\sca_\beta\cap W_\beta$.

We now apply (\ref{einstein}) to $E\in\End(\sg)$ defined by
$$
E:=\left[\begin{smallmatrix} 0&0\\ 0&\beta+||\beta||^2I\end{smallmatrix}\right],
$$
that is, $E|_{\ag}=0$ and $E|_{\ngo}=\beta+||\beta||^2I$.  If $\{ A_1,...,A_m\}$ is
an orthonormal basis of $\ag$ then the right hand side of (\ref{einstein}) is given
by
$$
\begin{array}{l}
\unc\sum\limits_{ij}\la E[e_i,e_j]-[Ee_i,e_j]-[e_i,Ee_j],[e_i,e_j]\ra \\
+\unc\sum\limits_{rs}\la E[A_r,A_s],[A_r,A_s]\ra \\
 + \unm\sum\limits_{ri}\la E[A_r,e_i],[A_r,e_i]\ra -
\unm\sum\limits_{ri}\la [A_r,Ee_i],[A_r,e_i]\ra,
\end{array}
$$
which in turn equals
\begin{equation}\label{prueba2}
\begin{array}{l}
\unc\la\pi(\beta+||\beta||^2I)\mu,\mu\ra + \unc\sum\limits_{rs}\la(\beta+||\beta||^2I)[A_r,A_s],[A_r,A_s]\ra \\
 + \unm\sum\limits_{ri}\la (\beta\ad{A_r}-\ad{A_r}\beta)(e_i),\ad{A_r}(e_i)\ra.
\end{array}
\end{equation}

The first and second terms in (\ref{prueba2}) are $\geq 0$ by Lemma \ref{delta} and
Lemma \ref{betapos}, respectively, and the last one equals
$\unm\sum\limits_r\la[\beta,\ad{A_r}],\ad{A_r}\ra$, which is $\geq 0$ by
(\ref{adbeta}) since $\ad{A_r}|_{\ngo}\in\Der(\mu)$ for all $r$.

We therefore obtain from (\ref{einstein}) and (\ref{c}) that
\begin{equation}\label{prueba1}
-\tfrac{\tr{S(\ad{H})^2}}{\tr{S(\ad{H})}}\tr{E}+\tr{S(\ad{H})E}\geq 0.
\end{equation}

Recall that $\tr{\beta}=-1$ and so
\begin{equation}\label{prueba3}
\begin{array}{ll}
\tr{E^2} & =\tr(\beta^2+||\beta||^4I+2||\beta||^2\beta)=||\beta||^2(1+n||\beta||^2-2) \\ \\
& = ||\beta||^2(-1+n||\beta||^2)=||\beta||^2\tr{E}.
\end{array}
\end{equation}

On the other hand, we have that
\begin{equation}\label{prueba4}
\tr{S(\ad{H})E}=\tr{\ad{H}|_{\ngo}(\beta+||\beta||^2)}=||\beta||^2\tr{S(\ad{H})}
\end{equation}
by Lemma \ref{betaort}.  We now use (\ref{prueba1}), (\ref{prueba3}) and
(\ref{prueba4}) and obtain
$$
\tr{S(\ad{H})^2}\tr{E^2}\leq (\tr{S(\ad{H})E})^2,
$$
a `backwards' Cauchy-Schwartz inequality.  This turns all inequalities mentioned
after (\ref{prueba2}) into equalities, in particular the second term:
$$
\unc\sum_{rs}\la(\beta+||\beta||^2I)[A_r,A_s],[A_r,A_s]\ra=0.
$$
We therefore get that $[\ag,\ag]=0$ since $\beta+||\beta||^2I$ is positive
definite by Lemma \ref{betapos}.

It only remains to consider the case $\mu=0$.  Here we argue in the same way but
with $E$ chosen as $E|_{\ag}=0$ and $E|_{\ngo}=I$.  It then follows from
(\ref{einstein}) that
$$
\sum_{rs}||[A_r,A_s]||^2=-\tfrac{\tr{S(\ad{H})^2}}{\tr{S(\ad{H})}}n+\tr{S(\ad{H})}
=\tfrac{\tr{S(\ad{H})^2}}{\tr{S(\ad{H})}}\left(\tfrac{(\tr{S(\ad{H})})^2}{\tr{S(\ad{H})^2}}-n\right)\leq
0,
$$
and thus $[\ag,\ag]=0$.  This concludes the proof of the theorem.
\end{proof}


\begin{thebibliography}{MMM}

\bibitem[AW]{AznWls} {\sc R. Azencott, E. Wilson}, Homogeneous manifolds with negative curvature I,
{\it Trans. Amer. Math. Soc.} {\bf 215} (1976), 323-362.

\bibitem[B]{Brg2} {\sc M. Berger}, A panoramic view of Riemannian geometry, Springer-Verlag,
Berlin-Heidelberg, 2003.

\bibitem[Bs]{Bss} {\sc A. Besse}, Einstein manifolds, {\it Ergeb. Math.} {\bf 10} (1987), Springer-Verlag,
Berlin-Heidelberg.

\bibitem[BWZ]{BhmWngZll} {\sc C. B$\ddot{{\rm o}}$hm, M.Y. Wang, W. Ziller}, A variational approach for
compact homogeneous Einstein manifolds, {\it Geom. Funct. Anal.} {\bf 14} (2004),
681-733.

\bibitem[D]{Dtt} {\sc I. Dotti}, Ricci curvature of left-invariant metrics on solvable unimodular Lie groups,
{\it Math. Z.} {\bf 180} (1982), 257-263.

\bibitem[FV]{Frd} {\sc H. Freudenthal, H. de Vries}, Linear Lie groups, {\it
Academic  Press}, New York and London (1969).

\bibitem[H]{Hbr} {\sc J. Heber}, Noncompact homogeneous Einstein spaces, {\it Invent. math}. {\bf 133} (1998), 279-352.

\bibitem[He]{Hss} {\sc W. Hesselink}, Uniform instability in reductive groups, {\it J. Reine Angew. Math.} {\bf 304}
(1978), 74-96.

\bibitem[Ke]{Kmp} {\sc G. Kempf}, Instability in invariant theory,
{\it Ann. Math.} {\bf 108} (1978), 299-316.

\bibitem[K]{Krw1} {\sc F. Kirwan}, Cohomology of quotients in symplectic and algebraic
geometry, {\it Mathematical Notes} {\bf 31} (1984), Princeton Univ. Press,
Princeton.

\bibitem[L1]{strata} {\sc J. Lauret}, On the moment map for the variety of Lie algebras,
{\it J. Funct. Anal.} {\bf 202} (2003), 392-423.

\bibitem[L2]{minimal} \bysame, A canonical compatible metric for geometric structures on
nilmanifolds, {\it Ann. Global Anal. Geom.} {\bf 30} (2006), 107-138.

\bibitem[L3]{cruzchica}  \bysame, Einstein solvmanifolds and nilsolitons, {\it Contemp. Math.} {\bf 491} (2009), 1-35.

\bibitem[LW]{einsteinsolv}  {\sc J. Lauret, C.E. Will}, Einstein solvmanifolds: existence
and non-existence questions, preprint 2006, {\it arXiv:} math.DG/0602502.

\bibitem[LeWa]{LbrWng} {\sc C. LeBrun, M.Y. Wang} (eds.), Surveys in differential geometry: essays on Einstein manifolds,
International Press, Boston (1999).

\bibitem[MFK]{Mmf} {\sc D. Mumford, J. Fogarty, F. Kirwan}, Geometric invariant theory, Third Edition,
Springer Verlag (1994).

\bibitem[N]{Nss} {\sc L. Ness}, A stratification of the null cone via the momentum
map, {\it Amer. J. Math.} {\bf 106} (1984), 1281-1329 (with an appendix by D. Mumford).

\bibitem[NN]{NktNkn} {\sc E.V. Nikitenko, Yu.G. Nikonorov}, Six-dimensional Einstein
solvmanifolds, {\it Siberian Adv. Math.} {\bf 16} (2006), 66-112.

\bibitem[Ni1]{Nkl} {\sc Y. Nikolayevsky}, Nilradicals of Einstein solvmanifolds, preprint 2006,
{\it arXiv:} math.DG/0612117.

\bibitem[Ni2]{Nkl3} \bysame, Einstein solvmanifolds and the
pre-Einstein derivation, {\it Trans. Amer. Math. Soc.}, in press ({\it arXiv:} math.DG/0802.2137).

\bibitem[RS]{RchSld} {\sc R.W. Richardson, P.J. Slodowy}, Minimum vectors for real reductive algebraic groups,
{\it J. London Math. Soc. (2)} {\bf 42} (1990), 409-429.

\bibitem[S]{Sch} {\sc D. Schueth}, On the standard condition for noncompact homogeneous Einstein spaces,
{\it Geom. Ded.} {\bf 105} (2004), 77-83.

\end{thebibliography}
\end{document}